\documentstyle[amssymb,amstex,graphicx,graphics,12pt,a4wide]{article}
\newtheorem{theorem}{Theorem}

\newtheorem{proposition}{Proposition}
\newtheorem{corollary}{Corollary}

\newcommand{\MaxFacets}{\operatorname{MaxPolytopes}_{d-1}}
\newcommand{\Cells}{\operatorname{Cells}}

\newcommand{\Vol}{\operatorname{Vol}}

\newcommand{\const}{\operatorname{const}}
\newcommand{\Var}{\operatorname{Var}}

\newcommand{\PLT}{\operatorname{PLT}}
\newcommand{\PHT}{\operatorname{PHT}}

\begin{document}

\title{\textbf{Second-Order Theory for Iteration Stable Tessellations}}
\author{by Tomasz Schreiber and Christoph Th\"ale\\ from Toru\'n and Fribourg}
\date{}
\maketitle

\begin{abstract}
This paper deals with iteration stable (STIT) tessellations, and, more generally, with a certain class of tessellations that are infinitely divisible with respect to iteration. They form a new, rich and flexible class of spatio-temporal models considered in stochastic geometry. The martingale tools developed in \cite{STP1} are used to study second-order properties of STIT tessellations. Firstly, a general formula for the variance of the total surface area of cell boundaries inside a convex observation window is shown. This general expression is combined with tools from integral geometry to derive explicit exact and asymptotic second-order formulae in the stationary and isotropic set-up, where a family of chord-power integrals plays an important role. Also a general formula for the pair-correlation function of the surface measure is found.
\end{abstract}
\begin{flushleft}\footnotesize
\end{flushleft}

\section{Introduction}

\begin{figure}[t]
\begin{center}
\includegraphics[width=7cm]{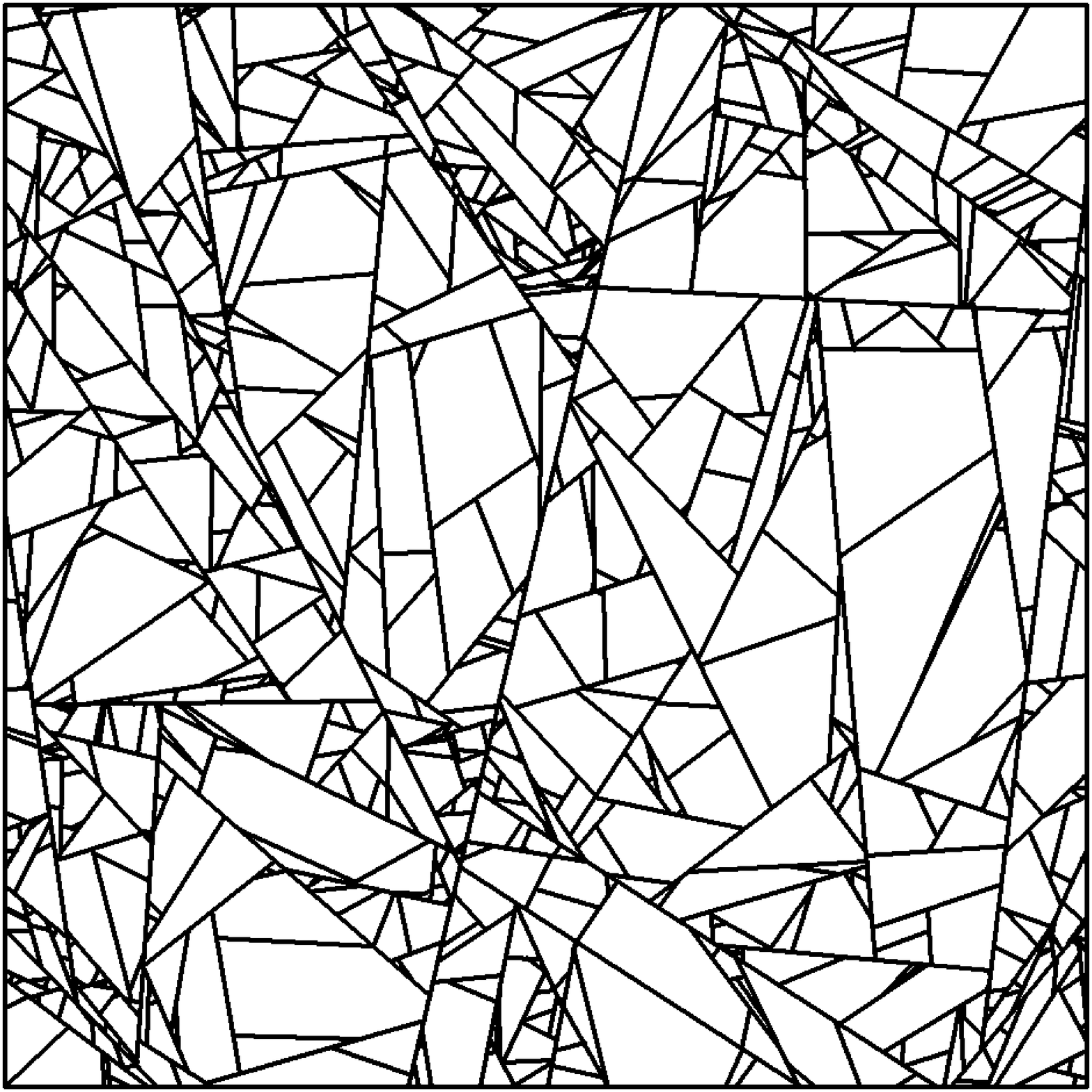}
\includegraphics[width=7cm]{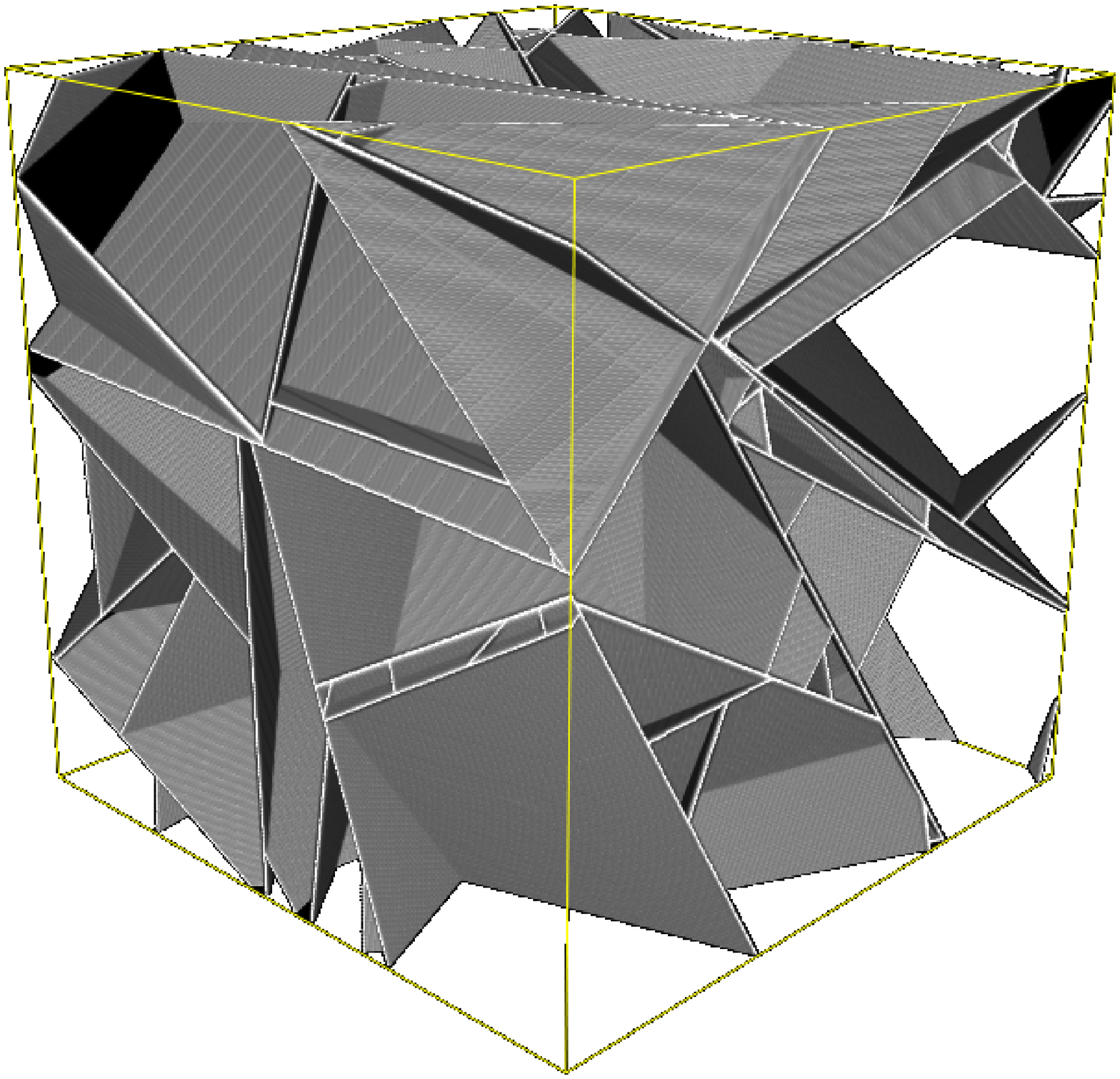}
\caption{Realizations of a planar and a spatial stationary and isotropic STIT tessellation (kindly provided by Joachim Ohser and Claudia Redenbach)}\label{Fig1}
\end{center}
\end{figure}

Iteration stable random tessellations (or mosaics), called STIT tessellations for short, form a new model for random tessellations of the $d$-dimensional Euclidean space and were formally introduced in \cite{MNW, MNW2, NW03, NW05}. They have quickly attracted considerable interest in stochastic geometry, because of their flexibility and analytical tractability. They clearly show the potential to become a new mathematical reference model beside hyperplane and Voronoi tessellations studied in classical stochastic geometry. Whereas much research in the last decades was devoted to mean values and mean value relations, modern stochastic geometry focusses on second-order theory and distributional results, see \cite{BL1, BL2, BR, HS0, HS1, HS2, HS3, S} to mention just a few.\\ To introduce the non-specialized reader to the subject, we briefly recall the basic construction of STIT tessellations within compact convex windows $W\subset{\Bbb R}^d$ with interior points. To this end, let us fix a (in some sense non-degenerate) translation-invariant measure $\Lambda$ on the space of hyperplanes. Further, let $t>0$ be fixed and assign to the window $W$ a random lifetime. Upon expiry of its lifetime, the primordial cell $W$ dies and splits into two sub-cells separated by a hyperplane hitting $W$, which is chosen according to the normalized distribution $\Lambda$. The resulting new cells are again assigned independent random lifetimes and the entire construction continues recursively until the deterministic time threshold $t$ is reached (see Figure \ref{Fig1} for an illustration). In order to ensure the Markov property of the above construction in the continuous-time parameter $t$, we assume from now on that the lifetimes are exponentially distributed. Moreover, we assume that the parameter of the exponentially distributed lifetimes of individual cells $[c]$ equals $\Lambda([c])$, where $[c]$ stands for the collection of hyperplanes hitting $c$. In this special situation, the random tessellation constructed by the described dynamics fulfils a stochastic \textbf{st}ability property under the operation of \textbf{it}eration of tessellations, and whence is indeed a \textit{STIT tessellations}. We refer to Section \ref{secSTIT} below for more details.\\ In \cite{STP1} we have introduced a new technique relying on martingale theory for studying these tessellations. One feature of this new approach is that it allows to investigate second-order parameters (i.e. variances) of the tessellation, which were out of reach so far and are in the focus of the present work. Based on a specialization of our martingale technique, we calculate in Section \ref{secSECOND} the variance of a general face-functional and as a special case we find the variance of the total surface area of cell boundaries in a bounded convex window. The resulting integral expression can be explicitly evaluated in the stationary and isotropic case by applying an integral-geometric transformation formula of Blaschke-Petkantschin type, which is also developed in this paper, see Section \ref{sdkfalsdkhf}. For the particular case of space dimension $3$, an exact formula without further integrals is found. Another important task in our context is to determine for fixed terminal times $t$ the large scale asymptotics of the afore calculated exact variance for a family of growing compact and convex windows with positive volume. Relying again on techniques from integral geometry, we will be able to determine asymptotic variance expressions, leading -- most interestingly -- in dimension $d=2$ to a result of very different qualitative nature compared with space dimensions $\geq 3$, where certain chord-power integrals, known from convex and integral geometry, will reflect the influence of the geometry of the observation window. In dimension $d=2$ we will see that, in contrast to the described situation for higher space dimensions, the shape of the window does not play any role and only its area enters our formulas, see Section \ref{ur9ngd8}. We also derive an explicit expression for the so-called pair-correlation function of the random surface measure for arbitrary space dimensions, see Section \ref{df8gd08fg7}, generalizing thereby recent findings from \cite{NOW}, which are based on completely different methods. This function is a a commonly used tool in spatial statistics and stochastic geometry to describe the second-order structure of a random set and it describes the expected surface density of the tessellation at a given distance from a typical point.\\ We would like to point out that the second-order theory developed in this paper is fundamental for our further work on STIT tessellations \cite{STP3, ST2, ST3} and that its extended version \cite{ST} is available online.\\ \\ In this paper we will make use of the following notation:
\begin{itemize}
    \item $B_R=B_R^d(o)$ is the $d$-dimensional ball around the origin with radius $R>0$.
    \item $\kappa_j := \Vol_j(B_1^j)$ is the volume of the $j$-dimensional unit ball, $j\kappa_j$ its surface area.	
    \item The uniform distribution on the unit sphere ${\cal S}_{d-1}$ in ${\Bbb R}^d$ (normalized spherical surface measure) is denoted by $\nu_{d-1}$.
\end{itemize}

\section{Construction and properties of the tessellations}\label{secSTIT}

Let $\Lambda$ be a non-atomic and locally finite measure on the space $\cal H$ of hyperplanes in the $d$-dimensional Euclidean space ${\Bbb R}^d$. Further, let $t>0$ and $W\subset{\Bbb R}^d$ be a compact convex window with interior points in which our construction of a random tessellation $Y(t\Lambda,W)$ is carried out. In a first step, we assign to the window $W$ an exponentially distributed random lifetime with parameter $\Lambda([W])$ where $[W] := \{ H \in {\cal H},\; H \cap W \neq \emptyset \}$ stands for the collection of hyperplanes hitting $W$. Upon expiry of its lifetime, the cell $W$ dies and splits into two polyhedral sub-cells $W^+$ and $W^-$ separated by a hyperplane in $[W]$, which is chosen according to the law $\Lambda(\cdot)/\Lambda([W]).$ The resulting new cells $W^+$ and $W^-$ are again assigned independent exponential lifetimes with respective parameters $\Lambda([W^+])$ and $\Lambda([W^-])$ (whence smaller cells live stochastically longer) and the entire construction continues recursively until the deterministic time threshold $t$ is reached (for an illustration see Figure \ref{Fig1}). The cell-separating $(d-1)$-dimensional facets (the word {\it facet} stands for a $(d-1)$-dimensional face here and throughout) arising in subsequent splits are usually referred to as $(d-1)$-dimensional maximal polytopes (or I-segments for $d=2$ as assuming shapes similar to the letter {\it I}).\\ The described process of recursive cell divisions is called the MNW-construction in honour of its inventors in the sequel and the resulting random tessellation created inside $W$ is denoted by $Y(t\Lambda,W)$ as mentioned above. The random tessellation $Y(t\Lambda,W)$ has the following properties (see \cite{NW05} for detailed proofs):
\begin{itemize}
 \item $Y(t\Lambda,W)$ is consistent in that $Y(t\Lambda,W) \cap V \overset{D}{=} Y(t\Lambda,V)$ for convex $V \subset W$ and thus $Y(t\Lambda,W)$ can be extended to random tessellation $Y(t\Lambda)$ on the whole space ${\Bbb R}^d$.
 \item If $\Lambda$ is translation-invariant, $Y(t\Lambda)$ is stationary, i.e. stochastically translation invariant. If, moreover, $\Lambda$ is the unit-density isometry-invariant hyperplane measure $\Lambda_{\rm iso}$, then $Y(t\Lambda_{\rm iso})$ is even isotropic, i.e. stochastically invariant under rotations wrt. the origin.
 \item $Y(t\Lambda,W)$ is iteration infinitely divisible with respect to the operation $\boxplus$ of \textbf{it}eration if tessellations for any compact convex $W\subset{\Bbb R}^d$. This is to say $$Y(t\Lambda,W)\overset{D}{=}m(Y((t/m)\Lambda,W)\boxplus\cdots\boxplus Y((t/m)\Lambda,W)),\ \ \ \ \ m=2,3,\ldots,$$ explained in detail in \cite{STP1}. Because of this property we call $Y(t\Lambda,W)$ an iteration infinitely divisible MNW-tessellation. In addition, if $\Lambda$ is translation-invariant, $Y(t\Lambda)$ is \textbf{st}able under the operation $\boxplus$, which is to say $$Y(t\Lambda)\overset{D}{=}m(Y((t/m)\Lambda)\boxplus\cdots\boxplus Y((t/m)\Lambda)),\ \ \ \ \ m=2,3,\ldots.$$ For this reason, $Y(t\Lambda)$ is called a random STIT tessellation in this case.
 \item In the stationary set-up, the surface density, i.e. the mean surface area of cell boundaries of $Y(t\Lambda)$ per unit volume equals $t$.
 \item STIT tessellations have the following scaling property: $$tY(t\Lambda)\overset{D}{=}Y(\Lambda),$$ i.e. the 
 tessellation $Y(t\Lambda)$ of surface intensity $t$ upon rescaling by factor $t$ has the same distribution as $Y(\Lambda)$, the STIT tessellation with surface intensity $1$.
\end{itemize}

\section{Background material}

In this section we recall a few facts from \cite{STP1}, wich are going to be crucial for our arguments below. Firstly, it follows directly from the MNW-construction of $Y(t\Lambda,W)$ that, in the continuous-time parameter $t$, this is a pure jump Markov process on the space of tessellations of $W$, whose generator is ${\Bbb L} := {\Bbb L}_{\Lambda; W}$ with
\begin{equation}\label{GEN}
 {\Bbb L}F(Y) = \int_{[W]} \sum_{f\in\Cells(Y\cap H)} [F(Y \cup \{f\}) - F(Y)] \Lambda(dH)
\end{equation}
for all $F$ bounded and measurable on space of tessellations of $W.$ Similar to the approch taken in \cite{STP1}, the general theory of Markov processes can now be used to construct a class of martingales associated with iteration infinitely divisible MNW-tessellations or, more specifically, STIT tessellations. Indeed, for bounded measurable $G = G(Y,t),$ considering the time-augmented Markov process $(Y(t\Lambda,W),t)_{t \geq 0}$ and applying standard theory, see Lemma 5.1 in Appendix 1 Sec. 5 in \cite{KipLan}, or simply by performing a direct check, we obtain
\begin{proposition}\label{propMART2} Assume that $G(Y,t)$ is twice continuously differentiable in $t$ and that $\sup_{Y,t} \left| \frac{\partial}{\partial t} G(Y,t) \right| + \left| \frac{\partial^2}{\partial t^2} G(Y,t) \right| < + \infty,$ which is condition (5.1) in \cite[App. 1 Sec. 5]{KipLan}. Then, the stochastic process
\begin{equation}\label{MART2}\nonumber
 G(Y(t\Lambda,W),t) - \int_0^t \left( [{\Bbb L}G(\cdot,s)](Y(s\Lambda,W)) + \frac{\partial}{\partial s}G(Y(s\Lambda,W),s) \right) ds
\end{equation}
is a martingale with respect to $\Im_t$, the filtration induced by $(Y(s\Lambda,W))_{0\leq s\leq t}$.
\end{proposition}
For $Y$ standing for some instant of $Y(t\Lambda,W)$, define
\begin{equation}\nonumber\label{FDEF}
 \Sigma_{\phi}(Y) := \sum_{f \in \MaxFacets(Y)} \phi(f)
\end{equation}
where, recall, $\MaxFacets(Y)$ are the $(d-1)$-dimensional maximal polytopes of $Y$ (the I-segments in the two-dimensional case), whereas $\phi(\cdot)$ is a generic bounded and measurable functional on $(d-1)$-dimensional facets in $W,$ that is to say a bounded and measurable function on the space of closed $(d-1)$-dimensional polytopes in $W,$ possibly chopped off by the boundary of $W,$ with the standard measurable structure inherited from space of closed sets in $W.$ Whereas the so-defined $\Sigma_\phi$ is not bounded, we cannot directly apply Dynkin's formula (see Appendix 1, Section 5 in \cite{KipLan} for example) to conclude that the stochastic process $\Sigma_\phi(Y(t\Lambda,W))-\int_0^t{\Bbb L}\Sigma_\phi(Y(s\Lambda,W))ds$ is a $\Im_t$-martingale. However, a suitable localization  argument can be applied (see \cite{STP1} for the details) to show this:
\begin{proposition}\label{propEXPECT} The stochastic process
$$\Sigma_{\phi}(Y(t\Lambda,W)) - \int_0^t \int_{[W]} \sum_{f\in\Cells(Y(s\Lambda,W)\cap H)} \phi(f) \Lambda(dH)ds$$
is a martingale with respect to $\Im_t$.
\end{proposition}

\section{A general variance formula}\label{secSECOND}

The general martingale statements from the previous section admit a convenient specialization to deal with second-order characteristics of iteration infinitely divisible MNW- or stationary STIT tessellations. Let us fix through this section a compact convex window $W\subset{\Bbb R}^d$ with interior points. From now on we will focus our attention on translation-invariant face functionals $\phi$ of $(d-1)$-dimensional facets, regarded as usual as closed subsets of $W$, of the form
\begin{equation}\label{PHIFORM}
 \phi(f) := \Vol_{d-1}(f) \zeta(\vec{\bf n}(f))
\end{equation}
with $\vec{\bf n}(f)$ standing for the unit normal to $f$ and $\zeta$ for a bounded measurable function on ${\cal S}_{d-1}.$ Recall now the definition of $\Sigma_\phi(Y(t\Lambda,W))$, introduce the bar notation $\bar\Sigma_{\phi}(Y(t\Lambda,W))$ by $\bar\Sigma_{\phi}(Y(t\Lambda,W)) := \Sigma_{\phi}(Y(t\Lambda,W)) - {\Bbb E}\Sigma_{\phi}(Y(t\Lambda,W))$ and put
\begin{equation}\label{APHI}
   A_{\phi}(Y(t\Lambda,W)) := \int_{[W]}\sum_{f\in\Cells(Y(t\Lambda,W)\cap H)}\phi(f)\Lambda(dH).
\end{equation}
Then we have
\begin{proposition}
The two stochastic process
\begin{equation}\label{MART3}
  \bar\Sigma_{\phi}(Y(t\Lambda,W)) \; \mbox{ and } \; \bar\Sigma_{\phi}^2(Y(t\Lambda,W)) - \int_0^t A_{\phi^2}(Y(s\Lambda,W)) ds
\end{equation}
are both $\Im_t$-martingales.
\end{proposition}
\paragraph{Proof of Proposition \ref{MART3}.}
For some instant $Y$ of $Y(t\Lambda,W)$ define $$G(Y,t) := (\Sigma_{\phi}(Y) - {\Bbb E}\Sigma_{\phi}(Y(t\Lambda,W)))^2,$$ so that $G(Y(t\Lambda,W),t) = \bar\Sigma_{\phi}^2(Y(t\Lambda,W))$. We use now Proposition \ref{propEXPECT} to check that
\begin{equation}\label{GPOCH}
 \frac{\partial}{\partial t} G(Y(t\Lambda,W),t) = - 2 [\Sigma_{\phi}(Y(t\Lambda,W)) - {\Bbb E}\Sigma_{\phi}(Y(t\Lambda,W))]
 {\Bbb E}A_{\phi}(Y(t\Lambda,W))
\end{equation}
with $A_\phi$ given by (\ref{APHI}). Put now together (\ref{GEN}), Proposition \ref{propEXPECT} and (\ref{GPOCH}) and use localization as in the discussion
preceding Proposition 2 in \cite{STP1} with $G_N,\; N \to \infty,$ chosen so that $(G_N(\cdot,\cdot) \wedge N) \vee -N \equiv (G(Y,t) \wedge N) \vee -N,$ that $|G_N(\cdot,\cdot)| \leq N+1$ and that $G_N(\cdot,t)$ be twice
continuously differentiable in $t,$ and with the localizing stopping times $T_N = \inf\{{t \geq 0}: 
(|G(Y(t\Lambda,W),t)| \vee |\frac{\partial}{\partial t}G(Y(t\Lambda,W),t)| \vee |\frac{\partial^2}{\partial t^2} G(Y(t\Lambda,W),t)|)
 \geq N\}.$ Proceeding as there, we readily conclude that
$$ \bar\Sigma_{\phi}^2(Y(t\Lambda,W)) - \int_0^t \int_{[W]} \sum_{f\in\Cells(Y(s\Lambda,W)\cap H)} \phi^2(f) \Lambda(dH) ds + $$
$$ 2 \int_0^t [ \int_{[W]} \sum_{f\in\Cells(Y(s\Lambda,W)\cap H)} \phi(f) [\Sigma_{\phi}(Y(s\Lambda,W))-{\Bbb E}\Sigma_{\phi}(Y(s\Lambda,W))] \Lambda(dH) - $$
$$ [\Sigma_{\phi}(Y(s\Lambda,W))-{\Bbb E}\Sigma_{\phi}(Y(s\Lambda,W))] {\Bbb E} A_{\phi}(Y(s\Lambda,W))] ds =  $$
\begin{equation}\label{VARMART}
 \bar\Sigma_{\phi}^2(Y(t\Lambda,W)) - \int_0^t A_{\phi^2}(Y(s\Lambda,W)) ds -
   2 \int_0^t \bar A_{\phi}(Y(s\Lambda,W)) \bar\Sigma_{\phi}(Y(s\Lambda,W)) ds 
\end{equation}
is a $\Im_t$-martingale, with $\bar A_{\phi}(Y(s\Lambda,W)) := A_{\phi}(Y(s\Lambda,W)) - {\Bbb E}A_{\phi}(Y(s\Lambda,W)).$ We now take advantage of the special form (\ref{PHIFORM}) of the face functional $\phi$ to conclude that $$ A_{\phi} \equiv \int_{[W]} \Vol_{d-1}(H \cap W) \zeta(\vec{\bf n}(H)) \Lambda(dH) = \const. $$ This implies $\bar A_{\phi} \equiv 0$ and thus, by Proposition \ref{propEXPECT} and (\ref{VARMART}), we can complete the proof.\hfill $\Box$\\ \\ The so-far established theory is now used to calculate the variance of face functionals as given by (\ref{PHIFORM}) of iteration infinitely divisible random MNW-tessellations $Y(t\Lambda,W)$ restricted to a compact convex window $W\subset{\Bbb R}^d$ with $\Vol_d(W)>0$.
\begin{theorem}\label{thmvariance} For arbitrary diffuse and locally finite measures $\Lambda$ on $\cal H$ and $\phi$ as in (\ref{PHIFORM}), we have $$\Var(\Sigma_{\phi}(Y(t\Lambda,W)))= \int_{[W]} \zeta^2(\vec{\bf n}(H)) \int_{W \cap H} \int_{W \cap H}\frac{1-\exp(-t\Lambda([xy]))}{\Lambda([xy])} dx dy \Lambda(dH).$$
\end{theorem}
\paragraph{Proof of Theorem \ref{thmvariance}.} Recall first (\ref{APHI}) and note that it implies
\begin{eqnarray}
\nonumber & & A_{\phi^2}(Y(t\Lambda,W))=\int_{[W]} \sum_{f \in \Cells(Y \cap H)} \phi^2(f) \Lambda(dH)\\
\nonumber &=&  \int_{[W]} \zeta^2(\vec{\bf n}(H)) \int_{W \cap H} \int_{W \cap H}
  {\bf 1}[x,y \mbox{ are in the same cell of } Y \cap H] dx dy 
  \Lambda(dH).\label{PHIKW}
\end{eqnarray} 
Thus, using (\ref{MART3}) and taking expectations of both sides yields immediately
 $$ \Var \Sigma_{\phi}(Y(t\Lambda,W)) =  \int_0^t \int_{[W]} 
    \zeta^2(\vec{\bf n}(H)) \int_{H \cap W} \int_{H \cap W} 
    {\Bbb P}(x,y \mbox{ are ...} $$
\begin{equation}\label{VAROG1}
  \mbox{... in the same cell of } Y(s\Lambda,W) \cap H) dx dy
    \Lambda(dH)ds. 
\end{equation}
Taking into account that $${\Bbb P}(x,y \mbox{ are in the same cell of } Y(s\Lambda,W) \cap H) = \exp(-s\Lambda([xy])),$$ which follows Theorem 1 \cite{STP1} together with standard properties of Poisson hyperplane tessellation, and using (\ref{VAROG1}), we end up with
$$
 \Var(\Sigma_{\phi}(Y(t\Lambda,W))) = \int_0^t \int_{[W]} 
    \zeta^2(\vec{\bf n}(H)) \int_{W \cap H} \int_{W \cap H} 
    \exp(-s\Lambda([xy])) dx dy \Lambda(dH) ds
$$
\begin{equation}\label{VAROG2}\nonumber  
= \int_{[W]} \zeta^2(\vec{\bf n}(H)) \int_{W \cap H} \int_{W \cap H} 
    \frac{1-\exp(-t\Lambda([xy]))}{\Lambda([xy])} dx dy \Lambda(dH),
\end{equation}
which completes our argument.\hfill $\Box$\\ \\ For general hyperplane measures $\Lambda$ this cannot be simplified further. However, in the special case, where $\Lambda$ is the unit-density isometry-invariant measure $\Lambda_{\rm iso}$, tools from integral geometry become available to evaluate the integral further.

\section{Exact variance expression for the isotropic STIT tessellation}\label{sdkfalsdkhf}

For the stationary and isotropic case $\Lambda = \Lambda_{\rm iso}$ we want to evaluate the variance expression from Theorem \ref{thmvariance} further in the special case $\phi=\Vol_{d-1}$, i.e. when $\zeta\equiv 1$. To simplify the notation we will write from now on $Y(t)$ instead of $Y(t\Lambda_{\rm iso})$.
\begin{theorem}\label{thmVAR} For the stationary and isotropic STIT tessellation $Y(t)$ with surface intensity $t>0$ we have 
\begin{equation}
\Var(\Vol_{d-1}(Y(t,W)))= \frac{d-1}{2} 
 \int_W \int_W {1-e^{-{2\kappa_{d-1}\over d\kappa_d}t\left\|x-y\right\|}\over\left\|x-y\right\|^2}dxdy  
\label{EQVAR0} \end{equation}
\begin{equation}
= {d(d-1)\kappa_d\over 2}\int_0^\infty\overline{\gamma}_W(r)r^{d-3}\left(1-e^{-{2\kappa_{d-1}\over d\kappa_d}tr}\right)dr,\label{EQVAR}
\end{equation}
where $W$ is a compact and convex subset of ${\Bbb R}^d$, and where $\overline{\gamma}_W(r)=\int_{{\cal S}_{d-1}}\Vol_d(W\cap(W+ru))\nu_{d-1}(du)$ is the isotropized set-covariance function of the window $W$.
\end{theorem}
The key to Theorem \ref{thmVAR} is a general integral-geometric transformation formula of Blaschke-Petkantschin type, which is interesting in its own right.
\begin{proposition}\label{propINTTRAFO} Let $W\subset{\Bbb R}^d$ be compact and convex and let $g:W\times W\rightarrow{\Bbb R}$ be a non-negative measurable function. Then
\begin{equation}
\int_{[W]}\int_{W\cap H}\int_{W\cap H}g(x,y)dxdy\Lambda_{\rm iso}(dH)={(d-1)\kappa_{d-1}\over d\kappa_d}\int_W\int_W{g(x,y)\over\left\|x-y\right\|}dxdy.\label{CALC1}
\end{equation}
\end{proposition}
\paragraph{Proof of Proposition \ref{propINTTRAFO}.} First, we use the affine Blaschke-Petkantschin formula \cite[Thm. 7.2.7]{SW} with $q=1$ to obtain for any non-negative measurable function $h:({\Bbb R}^d)^2\rightarrow{\Bbb R}$ $$\int_{{\Bbb R}^d}\int_{{\Bbb R}^d}h(x,y)dxdy={d\kappa_d\over 2}\int_{{\cal L}}\int_L\int_L h(x,y)\left\|x-y\right\|^{d-1}\ell_L(dx)\ell_L(dy)dL,$$ where ${\cal L}$ is the space of lines in ${\Bbb R}^d$ with invariant measure $dL$, i.e. the affine $1$-dimensional Grassmannian in ${\Bbb R}^d$ and $\ell_L$ is the the Lebesgue measure on $L$ with normalization as specified in \cite[Thm. 13.2.12]{SW}. Taking now $$h(x,y)={\bf 1}[x\in W]{\bf 1}[y\in W]\left\|x-y\right\|^kg(x,y)$$ for some $k>-d$ and another non-negative measurable function $g:W\times W\rightarrow{\Bbb R}^d$ we obtain
\begin{equation}
\int_W\int_W{\left\|x-y\right\|^kg(x,y)}dxdy={d\kappa_d\over 2}\int_{{\cal L}}\int_{W\cap L}\int_{W\cap L}\left\|x-y\right\|^{d-1+k}g(x,y)\ell_L(dx)d\ell_L(dy)dL.\label{EQBPF}
\end{equation}
For $k=-1$ this yields
\begin{equation}
\int_W\int_W{g(x,y)\over\left\|x-y\right\|}dxdy={d\kappa_d\over 2}\int_{{\cal L}}\int_{W\cap L}\int_{W\cap L}\left\|x-y\right\|^{d-2}g(x,y)\ell_L(dx)\ell_L(dy)dL.\label{CALC4}
\end{equation}
We replace now in (\ref{EQBPF}) for $k=0$, $W$ by $W\cap H$ for some fixed hyperplane $H$ and $d$ by $d-1$ and get
\begin{eqnarray}
\nonumber & & \int_{W\cap H}\int_{W\cap H}g(x,y)dxdy\\
\nonumber &=& {(d-1)\kappa_{d-1}\over 2}\int_{{\cal L}^H}\int_{W\cap H\cap L}\int_{W\cap H\cap L}\left\|x-y\right\|^{d-2}g(x,y)\ell_L(dx)\ell_L(dy)dL^H,
\end{eqnarray}
where by ${\cal L}^H$ we mean the $1$-dimensional affine Grassmannian restricted to $H$ with invariant measure $dL^H$ (this is the set of lines within hyperplane $H$). Averaging the last expression over all hyperplanes $H$ and using the fact that $\Lambda_{\rm iso}(dH)\otimes dL^H=dL,$ see \cite[Thm. 13.2.12]{SW}, yields \begin{eqnarray}
\nonumber & & \int_{{\cal H}}\int_{W\cap H}\int_{W\cap H}g(x,y)dxdy\Lambda_{\rm iso}(dH)\\
\nonumber &=& {(d-1)\kappa_{d-1}\over 2}\int_{{\cal H}}\int_{{\cal L}^H}\int_{W\cap H\cap L}\int_{W\cap H\cap L}\left\|x-y\right\|^{d-2}g(x,y)\ell_L(dx)\ell_L(dy)dL^H\Lambda_{\rm iso}(dH)\\
&=& {(d-1)\kappa_{d-1}\over 2}\int_{{\cal L}}\int_{W\cap L}\int_{W\cap L}\left\|x-y\right\|^{d-2}g(x,y)\ell_L(dx)\ell_L(dy)dL.\label{CALC3}
\end{eqnarray}
By comparing (\ref{CALC4}) and (\ref{CALC3}) we finally conclude
$$\int_{[W]}\int_{W\cap H}\int_{W\cap H}g(x,y)dxdy\Lambda_{\rm iso}(dH)={(d-1)\kappa_{d-1}\over d\kappa_d}\int_W\int_W{g(x,y)\over\left\|x-y\right\|}dxdy,$$
completing thereby the proof of the proposition.\hfill $\Box$
\paragraph{Proof of Theorem \ref{thmVAR}.}
In view of the general formula from Theorem \ref{thmvariance}, we take $$g(x,y)=\frac{1-\exp(-t\Lambda_{\rm iso}([xy]))}{\Lambda_{\rm iso}([xy])}={1-e^{-{2\kappa_{d-1}\over d\kappa_d}t\left\|x-y\right\|}\over {2\kappa_{d-1}\over d\kappa_d}\left\|x-y\right\|},$$ where the equality is a simple consequence of the mean projection formula from integral geometry, see \cite{SW}, Thm. 6.2.2 with $q=j=d-1$ there. Thus, upon applying the transformation formula (\ref{CALC1}) we conclude the following identity for $\Var(\Sigma_{\Vol_{d-1}}(Y(t,W)))=\Var(\Vol_{d-1}(Y(t,W)))$:
\begin{eqnarray}
\nonumber & & \int_{[W]}\int_{H\cap W}\int_{H\cap W}g(x,y)dxdy\Lambda_{\rm iso}(dH)={d-1\over 2}\int_W\int_W{1-e^{-{2\kappa_{d-1}\over d\kappa_d}t\left\|x-y\right\|}\over\left\|x-y\right\|^2}dxdy\\
\nonumber &=& {d(d-1)\kappa_d\over 2}\int_0^\infty\overline{\gamma}_W(r){1-e^{-{2\kappa_{d-1}\over d\kappa_d}tr}\over r^2}r^{d-1}dr\\
\nonumber &=& {d(d-1)\kappa_d\over 2}\int_0^\infty\overline{\gamma}_W(r)r^{d-3}\left(1-e^{-{2\kappa_{d-1}\over d\kappa_d}tr}\right)dr,
\end{eqnarray}
where we have passed to $d$-dimensional spherical coordinates.\hfill $\Box$\\ \\ In the special case $W=B_R^3$, the isotropized set-covariance function $\overline{\gamma}_{B_R^3}(r)$ takes the form $$\overline{\gamma}_{B_R^3}(r)=\begin{cases}{4\pi\over 3}R^3\left(1-{3r\over 4R}+{r^3\over 16R^3}\right) &: 0\leq r\leq 2R\\ 0 &: r>2R\end{cases}$$
and the variance integral can be evaluated in a closed form: \begin{equation}\Var(\Vol_{2}(Y(t,B_R^3)))={4\pi^2\over 3t^4}\left(t^2R^2(12-8tR+3t^2R^2)+24(1+tR)e^{-tR}-24\right).\label{EXPLICITVAR3D}\end{equation} The same closed form cannot be obtained for $d=2$, since $\overline{\gamma}_{B_R^2}(r)$ has a more complicated structure, i.e. $$\overline{\gamma}_{B_R^2}(r)=2R^2\arccos\left({r\over 2R}\right)-{r\over 2}\sqrt{4R^2-r^2}$$ for $r$ between $0$ and $2R$ and $\overline{\gamma}_{B_R^2}(r)=0$ for $r>2R$. Unfortunately, the resulting integral can in this case not further be simplified.

\section{The variance in the asymptotic regime}\label{ur9ngd8}

Another important task in our context is to determine for fixed $t$ the large $R$ asymptotics of the variance $\Var(\Vol_{d-1}(Y(t),W_R))$ for the family of growing windows $W_R=R\cdot W,\; R \to \infty$, with $W$ as in the previous section. Writing $\sim$ for the asymptotic equivalence of functions, i.e. $f(R)\sim g(R)$ iff $f(R)/g(R)\rightarrow 1$ as $R\rightarrow\infty$, we have
\begin{theorem}\label{CORASYFORMULA} For $d=2$,
\begin{equation}\label{EQVARASYM2D}
 \Var(\Vol_{1}(Y(t,W_R)) \sim \pi\Vol_2(W)R^2\log R,
\end{equation}
whereas for $d\geq 3$ we have
\begin{equation}\label{EQVARASYMHD}
\Var(\Vol_{d-1}(Y(t,W_R)))\sim R^{2(d-1)}{d-1\over 2}E_2(W)
\end{equation}
with $E_2(W)$ being the $2$-energy of $W$, see \cite[Chap. 8]{MATI}, given by
\begin{equation}\label{EN2}
 E_2(W) = \int_W \int_W \left\|x-y\right\|^{-2} dx dy. 
\end{equation}
\end{theorem}
In particular, this establishes weak long range dependencies present in stationary and isotropic STIT tessellations $Y(t)$. In the planar case, these dependencies are rather weak in that $${\Var(\Vol_{1}(Y(t,W_R)))\over\Vol_2(W_R)}\sim \pi\log R\rightarrow\infty,$$ as $R\rightarrow\infty$. For $d\geq 3$ these dependencies are much stronger, as the variance of the total surface area grows asymptotically like $R^{2(d-1)}$.
\paragraph{Proof of Theorem \ref{CORASYFORMULA}}
Formula (\ref{EQVARASYM2D}) can be established by using (\ref{EQVAR}), the relation $\overline{\gamma}_{W_R} \sim\Vol_2(W_R)=R^2\Vol_2(W)$,
valid uniformly for arguments $r = O(R / \log R),$ the observation that $\overline{\gamma}_{W_R} \to 0$ for $r = \Omega(R \log R),$ 
together with the fact that $\int_0^{L(R)}(1-e^{-cr}){dr\over r}\sim \log R$, $c>0$, as soon as $\log L(R) \sim \log R,$ and the scaling property of STIT tessellations: 
\begin{eqnarray}
\nonumber \Var(\Vol_1(Y(t,W_R))) &=& t^{-2}\Var(\Vol_1(Y(1,W_{tR}))) = {\pi\over t^2}\int_0^\infty\overline{\gamma}_{W_{tR}}(1-e^{-{2\over\pi}r}){dr\over r}\\
\nonumber &\sim& \pi t^{-2}\Vol_2(W_{tR})\log(tR)=\pi\Vol_2(W)R^2(\log R+\log t)\\
\nonumber &\sim& \pi R^2\Vol_2(W)\log R.
\end{eqnarray}
To see (\ref{EQVARASYMHD}), use (\ref{EQVAR0}) and again the scaling property of STIT tessellations to obtain
\begin{eqnarray}
\nonumber \Var(\Vol_{d-1}(Y(t,W_R))) &=& R^{2(d-1)} \Var(\Vol_{d-1}(Y(Rt,W))\\
\nonumber &=& R^{2(d-1)}\frac{d-1}{2} 
    \int_W \int_W {1-e^{-{2\kappa_{d-1}\over d\kappa_d}Rt\left\|x-y\right\|}\over\left\|x-y\right\|^2}dxdy\\
\nonumber &\rightarrow& R^{2(d-1)}{d-1\over 2}E_2(W),\ \ \ R\rightarrow\infty.
\end{eqnarray}
Observe that this does not extend for the separately treated case $d=2$ because there the integral in (\ref{EN2}) diverges.\hfill $\Box$\\ \\
It is easily seen that $E_2(\cdot)$ enjoys a \textit{superadditivity} property
\begin{equation}\label{SUPERADD}
\nonumber E_2(W_1 \cup W_2) \geq E_2(W_1) + E_2(W_2),\;\; W_1 \cap W_2 = \emptyset,
\end{equation}
which stands in contrast to (\ref{EQVARASYM2D}), where the asymptotic expression is linear in $\Vol_2(W).$
We will now derive an integral-geometric interpretation for this energy functional. Taking $g(x,y)\equiv 1$ and $k=-2$ in (\ref{EQBPF}) yields the identity
\begin{eqnarray}
\nonumber E_2(W) &=& \int_W\int_W{\left\|x-y\right\|^{-2}}dxdy={d\kappa_d\over 2}\int_{\cal L}\int_{W\cap L}\int_{W\cap L}\left\|x-y\right\|^{d-3}\ell_L(dx)\ell_L(dy)dL\\
 &=& {d\kappa_d\over(d-1)(d-2)}\int_{\cal L}\Vol_1(W\cap L)^{d-1}dL={2\over(d-1)(d-2)}I_{d-1}(W)\label{eee3comp}
\end{eqnarray}
with $I_{d-1}(W)$ being the $(d-1)$-st \textit{chord power integral} of $W$ in the sense of \cite[p. 363]{SW}. More precisely, $$I_{d-1}(W)={d\kappa_d\over 2}\int_{\cal L}\Vol_1^{d-1}(W\cap L)dL.$$ Hence, combining (\ref{eee3comp}) with (\ref{EQVARASYMHD}) from above ans using the fact that $I_{d-1}(\cdot)$ is homogeneous of degree $2(d-1)$, we arrive for $d\geq 3$ at
\begin{corollary} The asymptotic variance $\Var(\Vol_{d-1}(Y(t,W_R)))$, $R\rightarrow\infty$, is given by
\begin{equation}\nonumber\Var(\Vol_{d-1}(Y(t,W_R)))\sim {1\over d-2}I_{d-1}(W_R)={1\over d-2}R^{2(d-1)}I_{d-1}(W).\label{EQVARASYMHD2}\end{equation}
\end{corollary}
In general, $I_{d-1}(W)$ cannot further be evaluated. But for in the special case $W=B_1^d$ we have by applying \cite{SW}, Theorem 8.6.6 (with a corrected constant), $$I_{d-1}(B_1^d)={d2^{d-2}}{\kappa_d\kappa_{2d-2}\over\kappa_{d-1}}$$ and, thus, the 2-energy of the $d$-dimensional unit ball $B_1^d$ equals $$E_2(B_1^d)={d2^{d-1}\over (d-1)(d-2)}{\kappa_d\kappa_{2d-2}\over\kappa_{d-1}}={2\pi^d\over(d-1)(d-2)}\Gamma\left({d\over 2}\right)^{-2}.$$ In the particular case $d=3$ we obtain the value $E_2(B_1^3)=4\pi^2$, which agrees with our explicit variance formula (\ref{EXPLICITVAR3D}) from above.

\section{Pair-correlation function}\label{df8gd08fg7}
\begin{figure}[t]
\begin{center}
\includegraphics[width=8cm]{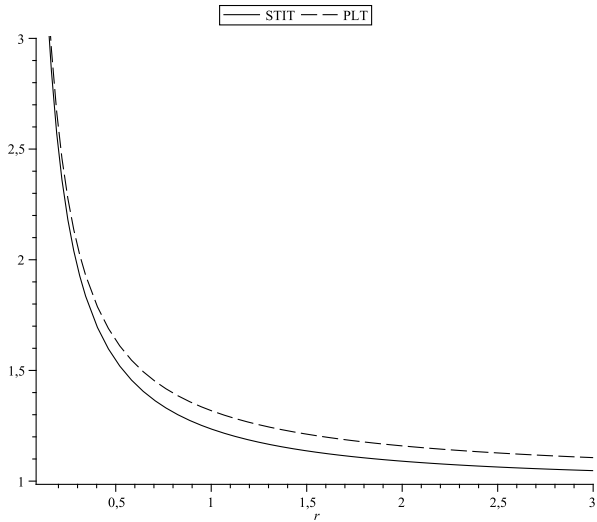}
\caption{Comparison of the pair-correlation function for a stationary and isotropic STIT tessellation $Y(1)$ and a stationary and isotropic Poisson line tessellation with edge length density $1$}\label{Fig2}
\end{center}
\end{figure}
It is our next goal to establish a closed formula for the pair-correlation function $g_d(r)$ of the random surface area measure of a STIT tessellation $Y(t)$, which is a commonly used tool in spatial statistics and stochastic geometry to describe the second-order structure of a random set. It describes the expected surface density of $Y(t)$ at a given distance $r$ from a typical point of $Y(t)$, see \cite{SW} or \cite{SKM} for exact definitions. It is well known (cf. \cite[p. 233]{SKM}) that the pair-correlation function $g_d(r)$ and the variance $\Var(\Vol_{d-1}(Y(t,W)))$ are related by the general formula $$\Var(\Vol_{d-1}(Y(t,W)))=d\kappa_dt^2\int_0^\infty\overline{\gamma}_W(r)(g_d(r)-1)dr.$$ Thus, from the explicit variance formula in Theorem \ref{thmVAR} the following can directly be deduced:
\begin{corollary}\label{corPCF} The pair-correlation function $g_d(r)$ of the random surface area measure of the stationary and isotropic random STIT tessellation $Y(t)$ is given by $$g_d(r)=1+{d-1\over 2t^2r^2}\left(1-e^{-{2\kappa_{d-1}\over d\kappa_d}tr}\right).$$
\end{corollary}
Especially for $d=2$, $g_d(r)$ becomes \begin{equation}\nonumber g_2(r)=1+{1\over 2t^2r^2}\left(1-e^{-{2\over\pi}rt}\right),\label{asdfkljasdf}\end{equation} which was independently obtained by Weiss, Ohser and Nagel by entirely different methods and is presented in \cite{NOW}. However, it should be emphasized though that our original approach developed above yields information also on higher dimensional cases. For example we have for the spatial case $d=3$, $$g_3(r) = 1+{1\over t^2r^2}\left(1-e^{-{1\over 2}tr}\right).$$ It is interesting to compare the pair-correlation function of $Y(t)$ from Corollary \ref{corPCF} with the corresponding function of a stationary and isotropic Poisson hyperplane tessellation $\PHT(t)$ with the same surface intensity $t>0$. The latter will be denoted by $g_d^{\PHT(t)}(r)$. Using Slivnyak's theorem for Poisson processes \cite[Thm. 3.3.5]{SW} one can easily show that $$g_d^{\PHT(t)}(r)=1+{(d-1)\kappa_{d-1}\over d\kappa_dtr}.$$ Especially for the planar case $d=2$, i.e. for the Poisson line tessellation abbreviated by $\PLT(t)$, we have $g_2^{\PLT(t)}(r)=1+1/(\pi tr).$ A comparison of $g_2(r)$ and $g_2^{\PLT(t)}(r)$ is shown in Figure \ref{Fig2}.

\subsection*{Acknowledgement}

The second author would like to thank Werner Nagel and Matthias Reitzner for comments and remarks.\\ The first author was supported by the Polish Minister of Science and Higher Education grant N N201 385234 (2008-2010). The second author was supported by the Swiss National Science Foundation grant SNF PP002-114715/1.

\begin{flushleft}
Tomasz Schreiber (1975-2010)\hfill Christoph Th\"ale\\
Faculty of Mathematics and Computer Science\hfill Department of Mathematics\\
Nicolaus Copernicus University\hfill University of Fribourg\\
Toru\'n, Poland\hfill Fribourg, Switzerland\\
\vspace{0.1cm}
\hfill \textit{Current address:}\\
\hfill Institute of Mathematics\\
\hfill University of Osnabr\"uck\\
\hfill Osnabr\"uck, Germany\\
\hfill christoph.thaele[at]uni-osnabrueck.de
\end{flushleft}
\end{document}